\documentclass{conm-p-l}
\usepackage{graphicx} % Required for inserting images
\usepackage{amssymb}
\usepackage[all,cmtip]{xy}
\input diagxy
\usepackage[backref]{hyperref}
\usepackage{todonotes}
\usepackage[alphabetic]{amsrefs}

\newtheorem{theorem}{Theorem}[section]
\newtheorem{proposition}[theorem]{Proposition}

\newtheorem{lemma}[theorem]{Lemma}
\theoremstyle{definition}
\newtheorem{definition}[theorem]{Definition}
\newtheorem{example}[theorem]{Example}

\newtheorem{remark}[theorem]{Remark}

\numberwithin{equation}{section}

\DeclareMathOperator{\Sk}{Sk}

\DeclareMathOperator{\Dec}{Dec}

\newcommand{\sm}{\mathsf{Smooth}}

\newcommand{\C}{\mathsf{C}}
\newcommand{\sSh}{\mathsf{sSh}}
\newcommand{\sC}{\mathsf{sC}}

\newcommand{\W}{W}
\newcommand{\LW}{\overline{W}}

\newcommand{\op}{op}

\DeclareMathOperator{\Cosk}{Cosk}

\DeclareMathOperator{\Tot}{Tot}
\DeclareMathOperator{\holim}{holim}

\newcommand{\MinnGp}{\mathsf{Min(n-Group)}}

\newcommand{\spKan}{\mathsf{Kan}^{split}_{\ast}}
\newcommand{\MinKan}{\mathsf{MinKan}_\ast}
\newcommand{\An}{\mathsf{An}}

\newcommand{\Id}{Id}

\newcommand{\Mfld}{\mathsf{Mlfds}}
\newcommand{\Schemes}{\mathsf{Schemes}}

\newcommand{\Sets}{\mathsf{Sets}}
\newcommand{\sSets}{\mathsf{sSets}}

\newcommand{\into}{\hookrightarrow}

\title{On Simplicial Principal Bundles in Descent Categories}
\dedicatory{For Ezra Getzler, on the occasion of his 60th birthday.}
\author{Jesse Wolfson}
\address{Department of Mathematics, University of California-Irvine}
\email{wolfson@uci.edu}

\thanks{The author was supported in part by NSF Grant DMS-1944862.}

\begin{document}

\maketitle
\begin{abstract}
Simplicial objects $\sC$ in descent categories $\C$, as introduced by Behrend and Getzler, provide a context in which to study higher stacks. In this note, we extend the construction of the canonical cocycle of a smooth principal $G$-bundle to the context of principal $G$-bundles in $\sC_{/X}$.  As an application, we  show how this specializes to $\C=\Sets$ to give a streamlined construction of $k$-invariants of reduced Kan complexes. We adapt this to give a similarly streamlined construction of minimal Kan complexes, with the goal of clarifying the role of the axiom of choice; Postnikov towers of minimal Kan complexes provide examples of towers of simplicial principal bundles of the type we consider.  
\end{abstract}

\section{Introduction}

Let $\C$ be a subcanonical descent category in the sense of \cite{BG}, or a ``category with covers'' in the sense of \cite{W}. Let $\sC$ denote the category of simplicial objects in $\C$.  Following \cite{H,P,W,BG,RZ}, we can view $\sC$ as a setting in which to study ``higher geometry'', i.e. higher Lie groupoids or geometric higher stacks in the geometry modeled by $\C$.  Geometric groups $G$ (Lie groups, algebraic groups, group schemes, etc.) and principal $G$-bundles form a basic tool in classical geometries, and (discrete) simplicial groups and their bundles have served a similar role in classical homotopy theory since their introduction in \cite{BGM}.  The goal of this short note is to flesh out a ``convex hull'' of these approaches.  Concretely, given a principal bundle
\[
    P\to X
\]
for a Lie group $G$ over a smooth manifold $X$, the {\em nerve} $NP$ of the bundle carries a canonical trivializing cocycle 
\[
    P\times_XP\cong P\times G.
\]
Via the implicit function theorem, it is enough to require that a principal bundle $P\to X$ in smooth manifolds be a surjective submersion in order to conclude that it is locally trivializeable in the open ball topology.  Our goal is to extend this package to $\sC$.  Actually, we fix $X\in \sC$ and work in the category $\sC_{/X}$ of maps $Y\to X$.  We refer to objects in this category as ``$X$-manifolds'', group objects $G\to X$ as ``$X$-groups'', etc.  In this setting, we take {\em covering left fibrations}, i.e. left fibrations
\[
    f\colon Z\to Y
\]
with $f_0\colon Z_0\to Y_0$ a cover, as the analogue of surjective submersions.\footnote{The ordinal category $\Delta$ comes equipped with an involution $t\colon \Delta\to\Delta$ given by ``reversing the orientation of simplices''.  Precomposing with this involution takes left fibrations into right fibrations and vice versa. For ease of exposition, we focus here on left fibrations.  We leave it as an exercise for the motivated reader to carry out the translation of this paper into right fibrations.} We now consider {\em Lie $X$-groups}, i.e. $X$-groups with $G\to X$ a level-wise cover.

\begin{example}
    Let $G$ be a Lie group acting by automorphisms on an abelian Lie group $A$.  Let $G$ act on itself by left translations and let $EG$ be the nerve of the associated action groupoid. Let $K(A,n)$ denote the simplicial abelian group whose normalized chain complex is $A$ in degree $n$ and 0 elsewhere.  Let $K(A,n)//G:=(K(A,n)\times EG)/G$ and consider the canonical map $K(A,n)//G\to EG/G=NG$. Then $K(A,n)//G\to NG$ is a Lie $NG$-group in the sense above.
\end{example}

\begin{definition}
    Let $Y\to X$ be an $X$-manifold, and $G\to X$ a Lie $X$-group. A {\em principal $G$-bundle} is a covering left fibration $P\to Y$ in $\sC_{/X}$ equipped with a $G$-action $\rho\colon G\times_X P\to P$ over $Y$ which gives an isomorphism of coequalizer diagrams
    \[
        \xymatrix{
             G\times_X P \ar@<-.5ex>[r]_{\pi_P} \ar@<.5ex>[r]^\rho \ar[d]^\cong_{\rho\times \pi_P} & P \ar@{=}[d] \ar[r] & Y \ar@{=}[d]\\
             P\times_Y P \ar@<-.5ex>[r]_{\pi_2} \ar@<.5ex>[r]^{\pi_1} & P \ar[r] & Y
        }.
    \]
\end{definition}

To define the nerve, we view both Lie $X$-groups $G\to X$ and the groupoid objects $Z\times_YZ\rightrightarrows Y$, for $Z\to Y$ a covering left fibration, as ``simplicial Lie groupoids''--i.e. groupoid objects in $\sC$ for which the source and target maps $d^h_1,d^h_0\colon \mathcal{G}_1\to \mathcal{G}_0$ are covers--and we denote these simplicial Lie groupoids by $G/X$ and $Z/Y$. Note that we do not require $\mathcal{G}_0$ to be a constant simplicial diagram, as is typical in the literature on simplicially enriched categories/groupoids. Following \cite{BGM}, we introduce functorial nerve $\LW$ and and universal bundle $\W$ constructions for simplicial Lie groupoids $\mathcal{G}$
\[
    \W\mathcal{G}\to\LW\mathcal{G}.
\]
For $p\colon Z\to Y$, the nerve comes with a natural augmentation
\begin{align*}
    \LW(p)\colon \LW(Z/Y)\to Y
\end{align*}

We can now state our main results.

\begin{proposition}\label{p:hyp}
    Let $f\colon Z\to Y$ be a covering left fibration.  The map
    \[
        \LW(f)\colon \LW(Z/Y)\to Y
    \]
    is a hypercover.
\end{proposition}

\begin{theorem}\label{t:can}
    Let $G\to X$ be a Lie $X$-group. Let $P\to Y$ be a principal $G$-bundle in $\sC_{/X}$. Then the nerve $\LW(P/Y)$ carries a canonical trivializing cocycle for $P$, i.e. there exists:
    \begin{enumerate}
        \item a natural map $\LW(P/Y)\to \LW(G/X)$, and
        \item a natural isomorphism 
        \[  
            \LW(P/Y)\times_Y P\to^\cong \LW(P/Y)\times_{\LW(G/X)} \W(G/X)
        \]
        over $\LW(P/Y)$
    \end{enumerate}
    where ``natural'' means, natural in $P\to Y$ and $G\to X$ in the usual sense.
\end{theorem}

\begin{remark}
    At the cost of introducing more formalism, we could equivalently restate Theorem~\ref{t:can} as the existence of a functor from principal $G$-bundles to a cocycle category in the sense of Jardine \cite{J}.
\end{remark}

In a final section, we specialize to the classical setting of Kan complexes, i.e. we take $\C=\Sets$ with surjections as covers, and consider the Moore-Duskin construction of the Postnikov tower of a reduced Kan complex $X$ as in \cite{H}. In this context, Proposition~\ref{p:hyp} and Theorem~\ref{t:can} give a streamlined construction of $k$-invariants, which we believe  should be adaptable to geometric contexts. In a similar spirit, we give a streamlined construction of minimal models of reduced Kan complexes (Theorem~\ref{t:min}) aimed at making more transparent the dependency of this construction on the axiom of choice. Reduced minimal Kan complexes provide natural examples of objects built from towers of principal bundles of the type we consider in Theorem~\ref{t:can}.

\subsection*{Acknowledgements}
The author thanks Ezra Getzler for helpful comments on a draft, and, beyond this, for introducing him to this circle of ideas, and for innumerable discussions, encouragement and advice over the years.  The author thanks the anonymous referee for helpful comments and corrections. This paper arose as an offshoot of an ongoing collaboration with Chris Rogers, and the author thanks him for many helpful conversations and advice.

\section{Preliminaries}\label{s:pre}
Except where otherwise noted, we follow the conventions of \cite{W}. Recall from \cite{BG} that a {\em subcanonical descent category} is a category $\C$ with a sub-category of ``covers'' such that
\begin{enumerate}
    \item $\C$ is closed under finite limits,
    \item pullbacks of covers are covers,
    \item if $f$ and $gf$ are covers, then so is $g$,
    \item covers are effective epimorphisms, i.e. if $f\colon X\to Y$ is a cover, then $X\times_Y X\rightrightarrows X\to^f Y$ is a coequalizer in $\C$.
\end{enumerate}
Examples of subcanonical descent categories include:
\begin{itemize}
    \item $\Sets$ with surjections as covers,
    \item $\An_{\mathbb{C}}$, the category of complex analytic varieties, with surjective submersions as covers,
    \item $\Schemes_S$, the category of schemes over a base $S$, with Zariski open covers (or surjective \'etale, fppf, or fpqc morphisms) as covers.
\end{itemize}
Subcanonical descent categories were introduced as a framework in which to study ``smooth'' maps (the covers) between possibly singular objects. From this perspective, the ``smooth'' objects are the objects $X\in\C$ such that the map $X\to \ast$ (from $X$ to the terminal object $\ast$) is a cover.  

If we wish to restrict our attention purely to ``smooth'' objects, we arrive at the framework of \cite{W}, which we informally refer to here as a ``category with covers.'' Concretely, this consists of a category $\C$ with a sub-category of ``covers'' subject to the following axioms:
\begin{enumerate}
    \item $\C$ has a terminal object $\ast$, and the map $X\to \ast$ is a cover for all $X\in\C$,
    \item pullbacks of covers along arbitrary maps exist and are covers,
    \item if $f$ and $gf$ are covers, then so is $g$,
    \item covers are effective epimorphisms.
\end{enumerate}
Examples include:
\begin{itemize}
    \item For any subcanonical descent category $\C$, the full subcategory of objects $X\in\C$ such that $X\to\ast$ is a cover, along with its subcategory of covers,
    \item $\Mfld$, the category of smooth Banach manifolds, with surjective submersions as covers,
    \item $\sm_k$, the category of smooth schemes over a field $k$, with surjective smooth morphisms as covers.
\end{itemize}
As principal bundles and torsors are important tools for studying both smooth and singular objects, our goal in the present note is to treat both contexts.  As a matter of exposition, we will take as a standing assumption that $\C$ is a category with covers--the geometric constructions one can make are more restricted in this setting--and we will comment as necessary on any adjustments for the case of a subcanonical descent category.

As in \cite{W}, the axioms for covers guarantee that they define a Grothendieck topology, and our assumption that covers are effective epimorphisms guarantees that the Yoneda embedding takes values in sheaves. We will freely identify an object in $\C$ with its functor of points, i.e. with the sheaf it represents, and we will talk about sheaves and (generalized) elements without further comment.  We write $\sC$ for the category of simplicial objects in $\C$, i.e. contravariant functors
\[
    X\colon \Delta^{\op}\to \C
\]
from the category of finite non-empty ordinals $\Delta$ into $\C$ and natural transformations between them, i.e. simplicial maps. We will identify objects in $\C$ as constant simplicial objects without further comment (i.e. we treat $\C$ as a full subcategory of $\sC$). For $X,Y\in\sC$, we write $\hom_{\sC}(X,Y)$ to denote the hom-sheaf, i.e. the sheaf on $\C$ whose value at $U\in \C$ is the set $\hom_{\sC}(U\times X,Y)$ of simplicial maps $U\times X\to Y$.\footnote{Note that we will not consider simplicial hom-sheaves, i.e. $\hom_{\sC}(X,Y)$ is a sheaf on $\C$, not a simplicial sheaf.} We view bisimplicial objects as ``horizontal'' simplicial diagrams in $\sC$ (where the internal simplicial direction in $\sC$ is ``vertical''). We write $[n]$ for the ordered set $\{0<\ldots<n\}$, $\Delta^n$ for the standard $n$-simplex, $\Lambda^n_i\subset \Delta^n$ for its ``$i^{th}$ horn'' and $\partial\Delta^n\subset\Delta^n$ for its boundary.  For $X\in\sC$, we write\footnote{Note that the right hand side involves a mild abuse of notation if coproducts of the terminal object $\ast$ do not define a fully faithful embedding $\Sets\to \C$. In this case, we view simplicial sets as simplicial objects in sheaves on $\C$ and use $\hom_{\sC}$ to denote $\hom_{\sSh(\C)}$.}
\begin{align*}
    \Lambda^n_iX&:=\hom_{\sC}(\Lambda^n_i,X)\\
    M_n X&:=\hom_{\sC}(\partial\Delta^n,X).
\end{align*}

For a map $f\colon Z\to Y$, we write $\Lambda^n_i(f)$ for the sheaf of commuting squares\footnote{N.b. equivalently, this can be described as the pullback sheaf $\Lambda^n_i Z\times_{\Lambda^n_i Y}Y_n$.}
\[
  \xymatrix{
    \Lambda^n_i \ar[r] \ar[d] & Z\ar[d] \\
    \Delta^n \ar[r] & Y
  }
\]
and $\lambda^n_i(f)\colon Z_n\to \Lambda^n_i(f)$ for the map induced by the inclusion $\Lambda^n_i\subset \Delta^n$.
\begin{definition}
    A map $f\colon Z\to Y$ in $\sC$ is a {\em left fibration} if for all $n$ and $0\le i<n$, the map
    \[
        \lambda^n_i(f)\colon Z_n\to \Lambda^n_i(f)
    \]
    is a cover.  We say it is a {\em covering left fibration} if, in addition, the map $f_0\colon Z_0\to Y_0$ is a cover.
\end{definition}

\begin{remark}
    Recall that a Kan fibration is a left fibration $f\colon Z\to Y$ such that the map
    \[
        \lambda^n_n(f)\colon Z_n\to \Lambda^n_n(f)
    \]
    is a cover for all $n$.  It is a {\em covering Kan fibration} if $f_0\colon Z_0\to Y_0$ is a cover.
\end{remark}

Following \cite{W,BG}, we develop some useful properties of left fibrations in $\C$. 
\begin{definition}\label{d:exp}
    Let $m\ge 0$.  A {\em left $m$-expansion} $S\into T$ be a map of simplicial sets for which there exists a filtration 
    \[
        S=F_{-1}T\subset F_0T\subset\cdots\subset F_r T=T
    \]
    such that there exists a sequence $n_\ell\ge m$, $0\le \ell\le r$, and a sequence $0\le i_\ell<n_\ell$ and pushout diagrams
    \[
        \xymatrix{
            \Lambda^{n_\ell}_{i_\ell} \ar[r] \ar[d] & F_{\ell-1}T\ar[d] \\
            \Delta^{n_\ell} \ar[r] & F_\ell T
        }.
    \]
    By a {\em left expansion} we mean a left 0-expansion.
\end{definition}

\begin{remark}
    Note that \cite{BG} imposes the additional requirement that the simplices in the expansion should monotonically increasing in dimension. While this does not change the theory, it leads, at times, to more combinatorially involved proofs (e.g. for Lemma~\ref{l:mexpcat} below).  We follow Whitehead \cite{Wh} and omit this requirement here.
\end{remark}

Immediately from Definition~\ref{d:exp}, we have the following.
\begin{lemma}\label{l:mexpcat}
    Let $R\into S$ be a left  $m_1$-extension and let $S\into T$ be a left $m_2$-extension.  Let $m=\min\{m_1,m_2\}$.  Then $R\into T$ is a left $m$-extension.
\end{lemma}

We now state a mild generalization of \cite[Lemma 3.8]{BG}.
\begin{lemma}\label{l:mexphorn}
    If $S$ is a union of $m\le n$ faces of the $n$-simplex $\Delta^n$, and $\partial_n\Delta^n\subset S$, then $S\into \Delta^n$ is a left expansion.
\end{lemma} 
\begin{proof}
    The proof is by induction on $n$.  For $n=1$ the statement is clear. For the induction step, enumerate the faces of $\Delta^n$ not in $S$ as $\{i_0,\ldots,i_{n-m}\}$ where $0\le i_{n-m}<\cdots<i_0<n$.  Let 
    \[
        F_\ell \Delta^n:=S\cup \bigcup_{j\le \ell} \partial_{i_j}\Delta^n.
    \]
    We claim that $F_{\ell-1} \Delta^n\cap\partial_{i_\ell}\Delta^n\into \partial_{i_\ell}\Delta^n$ is a left expansion for all $\ell<n-m$. To see this, observe that, by the simplicial identities, for any $0\le j<k\le n$, 
    \[
        \partial_{k-1}\partial_j\Delta^n=\partial_j\Delta^n\cap \partial_k \Delta^n=\partial_j\partial_k\Delta^n.
    \]
    By our assumption that $i_\ell<n$ for all $\ell$, $F_{\ell-1}\Delta^n\cap\partial_{i_\ell}\Delta^n$ is a union of at most $n-1$ faces of the $(n-1)$-simplex $\partial_{i_\ell}\Delta^n$ which contains the $(n-1)^{st}$ face. By the inductive hypothesis, it's a left expansion.  We conclude that the inclusion
    \[
        S\into F_{n-m-1}\Delta^n=\Lambda^n_{\iota_{n-m}}
    \]
    is a left expansion. Defining $F_{n-m}\Delta^n=\Delta^n$, we see that $F_{n-m-1}\Delta^n\into F_{n-m}\Delta^n$ is the canonical left expansion
    \[
        \Lambda^n_{\iota_{n-m}}\into \Delta^n.
    \]
    We now conclude that the map $S\into \Delta^n$ factors as a composition of left expansions, and thus by Lemma~\ref{l:mexpcat}, $S\into \Delta^n$ is a left expansion.
\end{proof}

\begin{lemma}\label{l:expcov}
    Let $S\into T$ be a left expansion.  Let $f\colon Z\to Y$ be a left fibration in $\sC$ such that the limit
    \[
        \hom(S,Z)\times_{\hom(S,Y)}\hom(T,Y)
    \]
    exists in $\C$. If $S\neq \emptyset$, then the limit $\hom(T,Z)$ exists and the map
    \[
        \hom(T,Z)\to \hom(S,Z)\times_{\hom(S,Y)}\hom(T,Y)
    \]
    is a cover. If $S=\emptyset$, then the same conclusion holds provided that $f$ is a covering left fibration.
\end{lemma}
\begin{proof}
    {\em Mutatis mutandis}, the proof is identical to the proof of \cite[Lemma 2.16]{W}.\footnote{Note that in {\em loc. cit.}, the edge case $S=\emptyset$ was mistakenly omitted from the proof.  Treating this edge case requires the hypothesis that $f_0$ is a cover in order to start the induction. }
\end{proof}

\begin{lemma}\label{l:covkan}
    Let $f\colon Z\to Y$ be a covering left fibration. Then for all $n$, the map $f_n\colon Z_n\to Y_n$ is a cover.
\end{lemma}
\begin{proof}
    We induct on $n$.  For $n=0$, this follows by assumption.  For the inductive step, for $f$ a covering left fibration, the map $f_n\colon Z_n\to Y_n$ factors as 
    \[
        Z_n\to Z_{n-1}\times_{Y_{n-1}}^{f,d_n}Y_n\to Y_n;
    \]
    %For $f$ a covering right fibration, we factor $f_n$ as 
    %\[
    %    Z_n\to Z_{n-1}\times_{Y_{n-1}}^{f,d_0}Y_n\to Y_n.
    %\]
    The second map is a pullback of the map $f_{n-1}\colon Z_{n-1}\to Y_{n-1}$, so the pullback $Z_{n-1}\times_{Y_{n-1}}Y_n$ exists in $\C$ and the map second map is a cover by the inductive hypothesis.  We can now apply Lemma~\ref{l:expcov} to conclude that the first map is a cover.
\end{proof}

For a map $f\colon Z\to Y$, we write $M_n(f)$ for the {\em $n^{th}$ relative matching object}, i.e. the sheaf of commuting squares\footnote{N.b. equivalently, this can be described as the pullback sheaf $M_n Z\times_{M_n Y} Y_n$.}
\[
    \xymatrix{
        \partial\Delta^n\ar[r] \ar[d] & Z\ar[d] \\
        \Delta^n \ar[r] & Y
    }
\]
and we write $\mu_n(f)\colon Z_n\to M_n(f)$ for the map induced by the inclusion $\partial\Delta^n\subset \Delta^n$.

\begin{definition}
    A map $f\colon Z\to Y$ in $\sC$ is a {\em hypercover} if for all $n$, the map 
    \[
        \mu_n(f)\colon Z_n\to M_n(f)
    \]
    is a cover.
\end{definition}
Hypercovers play the role of  trivial fibrations in the internal homotopy theory of $\sC$ (cf. \cite{BG,RZ}). In particular, they are covering Kan fibrations \cite[Theorem 2.17]{W}.  

We refer the reader to \cite{H,W,BG,RZ} for additional discussions of Kan fibrations, hypercovers, skeleta/coskeleta, representability issues, etc. 

\section{Nerves of Simplicial Lie Groupoids}
For a simplicial group $G$, Barratt, Gugenheim and Moore \cite{BGM} introduced a functorial nerve and universal bundle 
\[
    \W G\to \LW G
\]
with good homotopical properties in both discrete and geometric settings (cf. \cite[Section 6]{W}). 

Concretely:
\begin{align*}
    \W G_n&=G_n\times\cdots\times G_0\\
    \LW G_n&=G_{n-1}\times\cdots \times G_0
\end{align*}
and the face and degeneracy maps are given by multiplying an entry to the face of an adjacent one and inserting degeneracies.\footnote{c.f. \cite[Chapter IV.21]{May}, \cite[Chapter V.4]{GJ} or \cite[Section 6]{W}. Note that our conventions here follow \cite{May,GJ}, while the conventions in \cite{W} are equivalent, but distinct.  To pass from one to the other, pre-compose $\Tot$ by the transpose on bisimplicial objects (i.e. ``reflect across the diagonal''), and then interchange $\Dec_1$ with $\Dec^1$ in the definition of $\W$.} From a different perspective, we have natural isomorphisms of functors
\begin{align*}
    \LW &\cong \Tot\circ N\\
    \W   & \cong \Tot\circ\Dec_1\circ N
\end{align*}
where $N$ denotes the ``level-wise'' nerve, $\Tot$ denotes the Artin-Mazur totalization or ``co-diagonal'' (see \cite{AM,CR}), and $\Dec_1$ denotes the ``level-wise'' initial d\'ecalage  \cite{I} corresponding to the endo-functor on the ordinal category $\Delta$
\[
    [n]\mapsto [0]\sqcup [n].
\]

Our goal here is to generalize these to {\em simplicial Lie groupoids}.
\begin{definition}
    A {\em simplicial Lie groupoid} in $\C$ is a groupoid object $\mathcal{G}$ in $\sC$ for which both the source and target maps from the ``morphisms'' $\mathcal{G}_1$ to the ``objects'' $\mathcal{G}_0$
    \[
        d^h_1,d^h_0\colon \mathcal{G}_1\to \mathcal{G}_0
    \]
    are level-wise covers.
\end{definition}

Note that this is more general than the usual ``simplicially enriched groupoids'', as we do not require $\mathcal{G}_0$ to be simplicially constant. 

Our two primary cases of interest are:
\begin{enumerate}
    \item Lie $X$-groups $G\to X$, where the objects are just $X$, and the morphisms are $G$.\footnote{In a descent category, i.e. if $\C$ is closed under finite limits, then we can drop the assumption that $G\to X$ is a cover in the definition of the nerve below, as this is just there to ensure representability.}
    \item groupoids coming from level-wise covers $Z\to Y$, where the objects are just $Z$, and the morphisms are $Z\times_Y Z\to Z$.
\end{enumerate}
We recall our convention that all objects in $\sC$ are viewed as ``vertical'' simplicial diagrams.

\begin{definition}
    Let $\mathcal{G}$ be a simplicial Lie groupoid. Define its {\em nerve} $\LW\mathcal{G}$ by 
    \[
        \LW\mathcal{G}:=\Tot N \mathcal{G},
    \]
    where we apply $N$ level-wise to the simplicial Lie groupoid $\mathcal{G}$.
    
    Similarly, define the {\em universal bundle} $\W\mathcal{G}$ by
    \[
        \W\mathcal{G}:=\Tot \circ \Dec_1 \circ N \mathcal{G},
    \]
    with the projection
    \[
        \W\mathcal{G}\to \LW\mathcal{G}
    \]
    given by applying $\Tot$ to the component of the natural transformation $\Dec_1\to \Id$.
\end{definition}

The definition above has convenient functorial properties, which we will use later on. However, for computations, a more explicit presentation is often useful. 
\begin{lemma}
    Let $\mathcal{G}$ be a simplicial Lie groupoid.  There exists a natural isomorphism
    \[
        \LW\mathcal{G}_n\cong \mathcal{G}_{0,n}\times^{d^v_0,d^h_1}_{\mathcal{G}_{0,n-1}}\mathcal{G}_{1,n-1}\times^{d^h_0d^v_0,d^h_1}_{\mathcal{G}_{0,n-2}}\cdots\times^{d^h_0d^v_0,d^h_1}_{\mathcal{G}_{0,0}} \mathcal{G}_{1,0},
    \]
    and under this natural isomorphism, the degeneracy and face maps are given by
    \begin{align*}
        s_0(g_{0,n},g_{1,n-1},\ldots,g_{1,0})&=(s^v_0g_{0,n},s^h_0g_{0,n},g_{1,n-1},\ldots,g_{1,0})\\
        s_i(g_{0,n},g_{1,n-1},\ldots,g_{1,0})&=(s^v_ig_{0,n},s^v_{i-1}g_{1,n-1},\ldots,s^v_0g_{1,n-i},s^h_0d^h_0 g_{1,n-i},g_{1,n-i-1},\ldots,g_{1,0})\\
        d_0(g_{0,n},g_{1,n-1},\ldots,g_{1,0})&=(d^h_0g_{1,n-1},g_{1,n-2},\ldots,g_{1,0})\\
        d_n(g_{0,n},g_{1,n-1},\ldots,g_{1,0})&=(d^v_ng_{0,n},d^v_{n-1}g_{1,n-1},\ldots,d^v_1 g_{1,1})\intertext{and, for $0<i<n$,}
        d_i(g_{0,n},g_{1,n-1},\ldots,g_{1,0})&=(d^v_ig_{0,n},\ldots,d^v_1g_{1,n-i+1},d^v_0g_{1,n-i}\cdot g_{1,n-i-1},\ldots,g_{1,0}).
    \end{align*}
\end{lemma}

\begin{lemma}
    Let $\mathcal{G}$ be a simplicial Lie groupoid. There exists a natural isomorphism
    \[
        \W\mathcal{G}_n\cong \mathcal{G}_{1,n}\times^{d^h_0d^v_0,d^h_1}_{\mathcal{G}_{0,n-1}}\mathcal{G}_{1,n-1}\times^{d^h_0d^v_0,d^h_1}_{\mathcal{G}_{0,n-2}}\cdots\times^{d^h_0d^v_0,d^h_1}_{\mathcal{G}_{0,0}} \mathcal{G}_{1,0},
    \]
    and under this natural isomorphism, the degeneracy and face maps are given by
    \begin{align*}
        s_i(g_{1,n},g_{1,n-1},\ldots,g_{1,0})&=(s^v_ig_{1,n},\ldots,s^v_0g_{1,n-i},s^h_0d^h_0 g_{1,n-i},g_{1,n-i-1},\ldots,g_{1,0})\\
        d_i(g_{1,n},g_{1,n-1},\ldots,g_{1,0})&=(d^v_ig_{1,n},\ldots,d^v_1g_{1,n-i+1},d^v_0g_{1,n-i}\cdot g_{1,n-i-1},\ldots,g_{1,0})\\
        d_n(g_{1,n},g_{1,n-1},\ldots,g_{1,0})&=(d^v_ng_{1,n},d^v_{n-1}g_{1,n-1},\ldots,d^v_1 g_{1,1}).
    \end{align*}
    Under this isomorphism, the map $\W\mathcal{G}\to \LW\mathcal{G}$ is given on $n$-simplices by
    \[
        (g_{1,n},\ldots,g_{1,0})\mapsto (d^h_0 g_{1,n},g_{1,n-1},\ldots,g_{1,0}).
    \]
\end{lemma}

\begin{remark}
    We see from the lemma that the representability of $\W\mathcal{G}$ and $\LW\mathcal{G}$ follows directly from our assumption that source and target maps $d^h_0,d^h_1\colon \mathcal{G}_1\to \mathcal{G}_0$ are level-wise covers.\footnote{In a descent category, the same definition will also work for arbitrary groupoid objects in $\sC$.}
\end{remark}

\begin{example}
    Let $f\colon G\to X$ be a Lie $X$-group, viewed as a simplicial Lie groupoid. Let $\sigma\colon X\to G$ denote the identity section.  We denote its nerve by $\LW(G/X)$. Then 
    \begin{align*}
        \LW(G/X)_n&\cong X_n\times^{d_0,f}_{X_{n-1}}G_{n-1}\times^{fd_0,f}_{X_{n-2}}\cdots\times^{fd_0,f}_{X_0} G_0\intertext{and the face and degeneracy maps are given by}
        d_0(x,g_{n-1},\ldots,g_0)&=(f(g_{n-1}),g_{n-2},\ldots,g_0)\\
        d_i(x,g_{n-1},\ldots,g_0)&=(d_ix,d_{i-1}g_{n-1},\ldots,d_1g_{n-i+1},d_0g_{n-i}\cdot g_{n-i-1},g_{n-i-2},\ldots,g_0)\\
        d_n(x,g_{n-1},\ldots,g_0)&=(d_nx,d_{n-1}g_{n-1},\ldots,d_1g_1)\\
        s_0(x,g_{n-1},\ldots,g_0)&=(s_0x,\sigma(x),g_{n-1},\ldots,g_0)\\
        s_i(x,g_{n-1},\ldots,g_0)&=(s_ix,s_{i-1}g_{n-1},\ldots,s_0g_{n-i},\sigma f(g_{n-i}),g_{n-i-1},\ldots,g_0). 
    \end{align*}
    Taking $X=\ast$, we recover the usual formulas for $\LW$ of a simplicial group.
\end{example}

\begin{example}\label{ex:cov}
    Let $f\colon Z\to Y$ be a covering left fibration. Then $Z\times_Y Z\rightrightarrows Z$ gives a simplicial Lie groupoid by Lemma~\ref{l:covkan}, and we denote its nerve by $\LW(Z/Y)$. Then
    \begin{align*}
        \LW(Z/Y)_n&\cong Z_n\times^{fd_0,f}_{Y_{n-1}} Z_{n-1}\times^{fd_0,f}_{Y_{n-2}}\cdots\times^{fd_0,f}_{Y_0} Z_0\intertext{and face and degeneracy maps given by}
            d_0(z_n,\ldots,z_0)&=(z_{n-1},\ldots,z_0)\\
			d_i(z_n,\ldots,z_0)&=(d_iz_n,\ldots,d_1 z_{n-i+1},z_{n-i-1},\ldots,z_0)\\
			s_i(z_n,\ldots,z_0)&=(s_iz_n,\ldots,s_0z_{n-i},z_{n-i},\ldots,z_0).
    \end{align*}
    Note that in this case, we have a natural isomorphism
    \[
        \LW(Z/Y)\cong \Tot(\Cosk_{0,*} (Z)\times_{\Cosk_{0,*}Y} Y)
    \]
where $\Cosk_{0,*}(-)$ is the ``horizontal 0-coskeleton'' (i.e. which applies the 0-coskeleton to all the horizontal simplicial objects in a bisimplicial object).
\end{example}

\section{Proofs}
\begin{proof}[Proof of Proposition~\ref{p:hyp}]
    We claim that, for all $n$, there exists a commuting square
	\[
		\xymatrix{
			\bar{W}(Z/Y)_n \ar[rr]^{\mu_n(\LW(f))} \ar[d]_\cong && M_n\bar{W}(f) \ar[d]^\cong \\
			Z_n\times^{fd_0,f}_{Y_{n-1}} \bar{W}(Z/Y)_{n-1} \ar[rr]^{\lambda^n_0(f)\times 1} && \Lambda^n_0 (f)\times_{M_{n-1} Y}^{f\partial d_0,f\partial} \bar{W}(Z/Y)_{n-1}.
		}
	\]
	The proposition immediately follows from our assumption that $f$ is a covering left fibration (n.b. for $n=0$, the bottom horizontal arrow just becomes $Z_0\to Y_0$.)
	
	The existence of the above square reduces to an exercise in combinatorics. To see this, note that we have a natural isomorphism
	\[
		\bar{W}(Z/Y)\cong \Tot(\Cosk_{0,*} (Z)\times_{\Cosk_{0,*}Y} Y).
	\]
	The functor $\Tot$ is right adjoint to Illusie's ``total decalage'' $\Dec$ \cite{I,CR}, while the functor $\Cosk_{0,*}$ is right adjoint to the ``horizontal 0-skeleton''.  Similarly, the functor which takes a simplicial object and views it as a horizontally constant bisimplicial object is right adjoint to the colimit preserving functor given on bisimplices by
	\[
		s\otimes 1\colon \Delta^p\otimes \Delta^q\to \Delta^0\otimes\Delta^q
	\]
	where $s$ is the total degeneracy map.
	
	These natural isomorphisms combine to give a natural isomorphism between $\bar{W}(Z/Y)_n$ and the sheaf of commuting squares\footnote{N.b. this sheaf is equivalently the pullback \[\hom(\Sk_{0,*}\Dec\Delta^n,Z)\times_{\hom(\Sk_{0,*}\Dec\Delta^n,Y)}\hom(s\otimes 1 \Dec\Delta^n,Y).\]}
	\[
			\xymatrix{
				\Sk_{0,*}\Dec\Delta^n \ar[r] \ar[d] & Z \ar[d] \\
				(s\otimes 1) \Dec\Delta^n \ar[r] & Y
				} 
	\]
	Similarly, the definition of the matching object gives a natural isomorphism between $M_n\bar{W}(f)$ and the sheaf of commuting squares
	\[
			\xymatrix{
				\Sk_{0,*}\Dec\partial \Delta^n \ar[r] \ar[d] & Z \ar[d] \\
				(s\otimes 1) \Dec\Delta^n \ar[r] & Y
			}
	\]
	and the map $\mu_n(\bar{W}(f))$ corresponds to pre-composition by the commuting square
	\[
			\xymatrix{
			\Sk_{0,*}\Dec\partial \Delta^n \ar[r] \ar[d] & \Sk_{0,*}\Dec \Delta^n \ar[d]\\
			(s\otimes 1) \Dec\Delta^n  \ar@{=}[r] & (s\otimes 1) \Dec\Delta^n 
			}
	\]
	An exercise in bisimplicial combinatorics gives an isomorphism of this commuting square with
	\[
		\xymatrix{
			\Delta^0\otimes \Lambda^n_0\coprod \Sk_{0,*}\Dec\Delta^{n-1} \ar[r] \ar[d] & \Delta^0\otimes \Delta^n\coprod \Sk_{0,*}\Dec \Delta^{n-1} \ar[d]\\
			(s\otimes 1) \Dec\Delta^n  \ar@{=}[r] & (s\otimes 1) \Dec\Delta^n 
		}
	\]
	The original claim, and thus the proposition, immediately follows.
\end{proof}

\begin{proof}[Proof of Theorem~\ref{t:can}]
    Let $f\colon G\to X$ be a Lie $X$-group, $Y\to X$ an $X$-manifold and $P\to Y$ a principal $G$-bundle. Let $\pi\colon P\to Y\to X$ be the composite.  By definition, we have 
    \[
        \LW(P/Y):=\Tot N(P/Y).
    \]
    By definition, for all $n$, $P_n\to Y_n$ is a principal $G_n$-bundle in $\C_{/X}$, with the action compatible with the face and degeneracy maps.  The usual canonical cocycle on $N(P_n/Y_n)$ gives a natural isomorphism
    \[
        P_n\times_{X_n} \overbrace{G_n\times_{X_n}\cdots\times_{X_n}G_n}^m\to^\cong \overbrace{P_n\times_{Y_n}\cdots\times_{Y_n}P_n}^{m+1}=N(P_n/Y_n)_m=:N(P/Y)_{m,n}
    \]
    given on elements by
    \[
        (p,g_1,\ldots,g_m)\mapsto (p,g_1^{-1}p,\ldots,(g_1\cdots g_m)^{-1}p).
    \]
    Under these isomorphisms, the horizontal face and degeneracy maps of $N(P/Y)$ become
    \begin{align*}
        d^h_0(p,g_1,\ldots,g_m)&=(g_1^{-1}p,g_2,\ldots,g_m)\\
        d^h_i(p,g_1,\ldots,g_m)&=(p,g_1,\ldots,g_i\cdot g_{i+1},\ldots,g_n)\\
        s^h_0(p,g_1,\ldots,g_m)&=(p,\sigma f(g_1),g_1,\ldots,g_n)\\
        s^h_i(p,g_1,\ldots,g_m)&=(p,g_1,\ldots,g_i,\sigma f(g_i),g_{i+1},\ldots,g_n).
    \end{align*}
    The naturality of these isomorphisms implies that vertical face and degeneracy maps are computed coordinate-wise.  Applying $\Tot$, a straightforward exercise in the formulas gives a natural isomorphism
    \begin{align*}
        \LW(P/Y)_n&\to^\cong P_n\times^{\pi d_0,f}_{X_{n-1}} G_{n-1}\times^{fd_0,f}_{X_{n-2}}\cdots \times^{fd_0,f}_{X_0} G_0,\\
        (p_n,\ldots,p_0)&\mapsto (p_n,g_{n-1},\ldots,g_0)\intertext{where, for each $i$, $g_i\in G_i$ is the unique element such that $g_ip_i=d_0p_{i+1}$ (i.e. the canonical cocycle $(d_0p_{i+1},p_i)$ to $(g_i,p_i)$).  Under this isomorphism, the face and degeneracy maps are given by}
        d_0(p,g_{n-1},\ldots,g_0)&=(g_{n-1}^{-1}\cdot d_0p,g_{n-2},\ldots,g_0)\\
        d_i(p,g_{n-1},\ldots,g_0)&=(d_ip,d_{i-1}g_{n-1},\ldots,d_1g_{n-i+1},d_0g_{n-i}\cdot g_{n-i-1},g_{n-i-2},\ldots,g_0)\\
        s_0(p,g_{n-1},\ldots,g_0)&=(s_0p,\sigma\pi(p),g_{n-1},\ldots,g_0)\\
        s_i(p,g_{n-1},\ldots,g_0)&=(s_ip,s_{i-1}g_{n-1},\ldots,s_0g_{n-i},\sigma f(g_{n-i}),g_{n-i-1},\ldots,g_0).
    \end{align*}
    We conclude that the assignment 
    \[
        (p,g_{n-1},\ldots,g_0)\mapsto (\pi(p),g_{n-1},\ldots,g_0)
    \]
    gives a natural map
    \[
        \LW(P/Y)\to \LW(G/X)
    \]
    as desired.

    We now consider $P\times_Y \LW(P/Y)$.  Applying the isomorphism above, we obtain a natural isomorphism
    \[
        (P\times_Y\LW(P/Y))_n\cong P_n\times_{Y_n} P_n\times^{\pi d_0,f}_{X_{n-1}}G_{n-1}\times^{fd_0,f}_{X_{n-2}}\cdots\times^{fd_0,f}_{X_0}G_0.
    \]
    Post-composing this with the inverse of the canonical cocycle for the $G_n$-principal bundle $P_n\to Y_n$
    \[
        (g_n,p)\mapsto (g_np,p),
    \]
    we obtain a natural isomorphism
    \begin{align*}
        (P\times_Y \LW(P/Y))_n&\to^\cong P_n\times_{X_n} G_n\times_{X_{n-1}} \cdots \times_{X_0} G_0\\
        (p,g_n^{-1}p,g_{n-1},\ldots,g_0)&\mapsto (p,g_n,g_{n-1},\ldots,g_0)
    \end{align*}
    which by inspection gives a simplicial isomorphism over $\LW(P/Y)$, as desired.
\end{proof}
\section{Applications}
We close by showing how Theorem~\ref{t:can} gives streamlined constructions of $k$-invariants and minimal models in the classical context of Kan simplicial sets. Let $\C=\Sets$ with surjections as covers.  For each $n\ge 0$, Moore \cite{M} and Duskin \cite{D} introduced ``truncation'' functors
\[
    \tau_{<n},\tau_{\le n}\colon \sSets\to\sSets
\]
such that for, $X$ a Kan complex, we have a natural tower of Kan complexes
\[
    X\to \cdots\to\tau_{\le n}X\to \tau_{<n}X\to^\sim \tau_{\le n-1}X\to\cdots\to \tau_{<1}X\to^\sim \tau_{\le 0}X=\pi_0X,
\]
with 
\[
    \tau_{<n+1}X\to^\sim \tau_{\le n}X
\]
a trivial fibration, and 
\[
    \tau_{\le n}X\to \tau_{<n}X
\]
a minimal Kan fibration for all $n$, and with
\[
    X=\varprojlim_n \tau_{\le n}X=\varprojlim_n\tau_{<n}X\simeq\holim_n \tau_{\le n}X\simeq \holim_n\tau_{<n}X.
\]
See \cite[Section II.8]{May} and \cite[Chapter IV.3]{GJ} for a discussion of $\tau_{<n}$ in simplicial sets.  See \cite{H,W} for a discussion of both functors in geometric contexts.

For $X$ reduced\footnote{We assume reduced for simplicity. The discussion goes though {\em mutatis mutandis} for $X$ connected with a chosen basepoint.}, by \cite[Proposition V.2.2]{BGM}, 
\[
    \tau_{\le n}X\to \tau_{<n}X
\]
is a minimal Kan fiber bundle, and in fact a twisted Cartesian product of the Eilenberg-MacLane space $K(\pi_n X,n):=\LW^n(\pi_n X)$ with twisting function given by the action of $\pi_1X$ (cf. \cite{BGM} or \cite[Chapter VI.5]{GJ}). Checking the definitions, we see this is in fact a principal bundle for the $N\pi_1X$-group 
\[
    K(\pi_n X,n)//\pi_1X:=(\W \pi_1 X\times K(\pi_n X,n))/\pi_1 X.
\]
Applying Theorem~\ref{t:can}, we obtain a commuting diagram in the category $\sSets_{/N\pi_1X}$ for each $n\ge 1$
\[
    \xymatrix{
        \tau_{\le n}X\times_{\tau_{<n}X}\LW(\tau_{\le n}X/\tau_{<n}X) \ar[r] \ar[d] & \W(K(\pi_n X,n)//\pi_1X)\ar[d] \\
        \LW(\tau_{\le n}X/\tau_{<n}X) \ar[r] \ar@{->>}[d]_\sim &  \LW(K(\pi_n X,n)//\pi_1X)\\
        \tau_{<n}X
    }
\]
where the bottom vertical map is a trivial fibration, by Proposition~\ref{p:hyp}. Using the canonical isomorphism 
\[
    \LW(K(\pi_nX,n)//\pi_1X)\cong K(\pi_n X,n+1)//\pi_1X
\]
the above diagram is naturally isomorphic to
\[
    \xymatrix{
        \tau_{\le n}X\times_{\tau_{<n}X}\LW(\tau_{\le n}X/\tau_{<n}X) \ar[r] \ar[d] & \W K(\pi_n X,n)//\pi_1X\ar[d] \\
        \LW(\tau_{\le n}X/\tau_{<n}X) \ar[r] \ar@{->>}[d]_\sim &  K(\pi_n X,n+1)//\pi_1X\\
        \tau_{<n}X
    },
\]
i.e. this is precisely an $n$-dimensional $k$-invariant for the Kan complex $X$ in the sense of \cite[Proposition VI.5.1]{GJ}. Comparing this to the construction in \cite[Chapter VI.5]{GJ}, we see that the above is essentially the construction of {\em loc. cit.}, but in a streamlined fashion which makes clear the categorical requirements. In particular, one only needs the representability of the functors $\tau_{\le n}$ and $\tau_{<n}$, e.g. as discussed in \cite{H,W}, and the analogue of the results of \cite[Sections III, IV]{BGM} to obtain a functorial theory of $k$-invariants for reduced Kan objects in geometric contexts $\sC$.

To put this in context, we now show the following. 
\begin{theorem}\label{t:min}
    Let $\spKan$ denote the category whose objects tuples $(X,\{\sigma_n\}_{n\ge 0})$ where $X$ is a reduced Kan complex, and for all $n\ge 0$, $\sigma_n\colon \tau_{\le n}X\to \tau_{<n+1}X$ is a section of the trivial fibration $\tau_{<n+1}X\to^\sim \tau_{\le n}X$.\footnote{Note that such sections always exist by the axiom of choice and the definition a trivial fibration of simplicial sets.}  Then there exists a functor
    \[
        \mu\colon \spKan\to \MinKan
    \]
    taking values in reduced {\em minimal} Kan complexes and a natural weak equivalence
    \[
        \alpha\colon \mu(X,\{\sigma_n\})\to^\sim X
    \]
    which is an inclusion of simplicial sets.
\end{theorem}
\begin{proof}
    Let $(X,\{\sigma_n\})\in \spKan$. Let 
    \[
        \mu_0(X,\{\sigma_n\}):=\pi_0X=\ast
    \]
    and let $\alpha_0\colon \mu_0(X,\{\sigma_n\})\to^1 \tau_{\le 0}X$ be the identity map. Now suppose that we have constructed a functor
    \[
        \mu_n\colon \spKan\to \MinnGp
    \]
    taking values in the category of minimal $n$-groups, i.e. the full sub-category $\MinKan$ consisting of those $X$ such that $X=\tau_{\le n}X$. Suppose also that we have constructed a natural transformation as above
    \[
        \alpha_n\colon \mu_n(X,\{\sigma_n\})\to \tau_{\le n} X,
    \]
    such that $\tau_{\le n-k}\alpha_n=\alpha_{n-k}$ for all $k\le n$. Then we define $\mu_{n+1}(X,\{\sigma_n\})$ to be the pullback
    \[
        \xymatrix{
            \mu_{n+1}(X,\{\sigma_n\}) \ar[rr] \ar[d] && \tau_{\le n+1}X \ar[d]\\
            \mu_n(X,\{(\sigma_n)\}) \ar[r]^{\alpha_n}_\sim & \tau_{\le n}X \ar[r]^{\sigma_n}_\sim & \tau_{<n+1}X
        }
    \]
    and we define $\alpha_{n+1}\colon \mu_{n+1}(X,\{\sigma_n\})\to \tau_{\le n+1}X$ to be the top horizontal map.  Because $\mu_n(X,\{(\sigma_n)\})$ is a minimal reduced Kan complex and the right vertical map is a minimal Kan fibration with reduced and $n$-connected fiber, we conclude that $\mu_{n+1}(X,\{\sigma_n\})$ is a minimal reduced Kan complex and the left vertical map is the $n$-truncation.  Because the right vertical map is a Kan fibration and the bottom horizontal map is a weak equivalence and an injection, we conclude that $\alpha_{n+1}$ is a weak equivalence and an injection.  By the universal property of pullbacks, the naturality of $\alpha_n$ and the definition of morphisms in $\spKan$, we conclude that $\alpha_{n+1}$ is natural transformation as claimed. This completes the induction step. 
    Defining
    \[
        \mu(X,\{\sigma_n\})=\varprojlim_n \mu_n(,\{(\sigma_n)\})
    \]
    and 
    \[
        \alpha=\varprojlim_n \alpha_n
    \]
    and noting, from the above square, that $\tau_{\le n}\mu_{n+1}=\mu_n$, we conclude that $\mu$ takes values in reduced minimal Kan complexes, that $\tau_{\le n}\mu=\mu_n$.  Noting that $X=\varprojlim_n \tau_{\le n}X$, we conclude that 
    \[
        \alpha\colon \mu(X,\{\sigma_n\})\to X
    \]
    is a natural weak equivalence and an injection as desired.
\end{proof}

\begin{remark}
The above presents a streamlined version of the classical construction of minimal models (cf. \cite[Chapter I.10]{GJ}) which we believe clarifies the role of the axiom of choice.  While it does not strictly require the data of the $k$-invariants constructed above, we include it here as an example of a tower of principal $G$-bundles for $Z$-groups $G\to Z$, to which Theorem~\ref{t:can} applies.
\end{remark}

\end{document}